%% file: FreeBF2013.tex
\documentclass[12pt,a4paper]{article}
\usepackage[T1]{fontenc}
\usepackage{lmodern}
\usepackage{amsmath,amssymb,amsfonts}        
\usepackage{mathptmx}
\usepackage[scaled=0.91]{helvet}
\usepackage{courier}
\usepackage[dvips]{graphicx}
\usepackage{psfrag}
\usepackage{color}
\usepackage{url}

\newtheorem{theorem}{Theorem}
\def\ds{\displaystyle}
\def\e{\epsilon}
\def\FBF{free boundary formulation}

\def\www{.9\textwidth}

\title{Free Boundary Formulation for
 BVPs on a Semi-Infinite Interval and Non-Iterative Transformation Methods}
\author{Riccardo Fazio \\
Department of Mathematics and Computer Science\\
University of Messina \\
Viale F. Stagno D'Alcontres, 31 \\
98166 Messina, Italy \\
E-mail: rfazio@unime.it \\
Home-page: \url{http://mat521.unime.it/fazio}}
\date{\today}
\begin{document}
\maketitle
\begin{abstract}
This paper is concerned with two examples on the application of the free boundary formulation to BVPs on a semi-infinite interval.
In both cases we are able to provide the exact solution of both the BVP and its \FBF .
Therefore, these problems can be used as benchmarks for the numerical methods applied to BVPs on a semi-infinite interval and to free BVPs. 
Moreover, we emphasize how for two classes of free BVPs, we can define non-iterative initial value methods, whereas BVPs are usually solved iteratively.
These non-iterative methods can be deduced within Lie's group invariance theory.
Then, we show how to apply the non-iterative methods to the two introduced \FBF s in order to obtain meaningful numerical results. 
Finally, we indicate several problems from the literature where our non-iterative transformation methods can be applied.
\end{abstract}
\smallskip

\smallskip



\section{Introduction}
Usually when dealing for the first time with a boundary value problem (BVP) defined on a semi-infinite interval the applied scientist does not know the exact or even an approximate solution.
As a consequence, he often is tempted to try for a numerical solution of the problem.
Therefore, along the years several approaches have been proposed  in order to solve BVPs defined on a semi-infinite interval numerically.

The oldest and simplest approach is to replace the original problem by one defined on a finite interval, where a finite value, the so-called truncated boundary, is used instead of infinity. 
This approach was used, for instance, by Horwarth \cite{Horwarth:1938:SLB} and by Goldstein \cite[p. 136]{Goldstein:1938:MDF} to get the tabulated numerical solution of the Blasius problem \cite{Blasius:1908:GFK}.
However, to get an accurate solution a comparison of numerical results obtained for several values of the truncated boundary is necessary as suggested by Fox \cite[p. 92]{Fox} or by Collatz \cite[pp. 150-151]{Collatz}.
Moreover, in some cases accurate solutions can be found only by using very large values of the truncated boundary.
This is, for instance, the case for the branches of the von Karman swirling flows where values of truncated boundaries up to several hundreds were used by Lentini and Keller \cite{Lentini:KSF:1980}.

The limitation of the above classical approach has lead de Hoog and Weiss \cite{deHoog:1980:ATB}, Lentini and Keller \cite{Lentini:BVP:1980} and Markowich
\cite{Markowich:TAS:1982,Markowich:ABV:1983} to develop a theory for defining the asymptotic boundary conditions to be imposed at a given value of the truncated boundary. 
Those asymptotic boundary conditions have to be derived by a preliminary asymptotic analysis involving the Jacobian matrix of  the right-hand side of the governing equations evaluated at infinity.
The main idea of this asymptotic boundary conditions approach is to project the solution into the manifold of bounded solutions.
By using the same value of the truncated boundary, a more accurate numerical solution can be found by this approach than that obtained by the classical approach, because in the first case the imposed boundary conditions are obtained from the asymptotic behaviour of the solution.
However, we should note that this approach is not straightforward, see the remarks in \cite{Ockendon}, and that for nonlinear problems highly nonlinear asymptotic boundary conditions usually result.
Asymptotic boundary conditions have been applied successfully to the numerical approximation of the so-called \lq \lq connecting orbits\rq \rq \ problems of dynamical systems, see Beyn \cite{Beyn:1990:GBN,Beyn:1990:NCC,Beyn:1992:NMD}.
Those problems are of interest, not only in connection with dynamical systems, but also in the study of travelling wave solutions of partial differential equations of parabolic and hyperbolic type as shown by Beyn \cite{Beyn:1990:NCC}, Friedman and Doedel \cite{Friedman:1991:NCC}, Bai et al. \cite{Bai:1993:NCH}, and Liu et al. \cite{Liu:1997:CCH}.

A different approach, for the numerical solution of BVPs defined on a semi-infinite, is to consider a free boundary formulation of the given problem, where the unknown free boundary can be identified with a truncated boundary. 
In this approach the free boundary is unknown and has to be found as part of the solution. 
This free boundary approach overcomes the need for a priori definition of the truncated boundary. 
Free BVPs represent a numerical challenge because they are always nonlinear as pointed out first by Landau \cite{Landau:1950:HCM}.
However, a \FBF \ has been successfully applied to several problems in the applied sciences: namely,  the Blasius problem \cite{Fazio:1992:BPF}, a two-dimensional stagnation point flow \cite{Ariel:FEF:1993},
the Falkner-Skan model \cite{Fazio:1994:FSEb,Zhang:IMS:2009,Zhu:NSF:2009}, and the model describing a fluid flowing around a slender parabola of revolution \cite{Fazio:1996:NAN} in boundary layer theory, the computation of a two-dimensional homoclinic connecting orbit \cite{Fazio:2002:SFB}, and a problem related to the deflection of a semi-infinite pile embedded in soft soil \cite{Fazio:2003:FBA}.
The last problem is of interest in foundation engineering, for instance, in the design of drilling rigs above the ocean floor, see Lentini and Keller \cite{Lentini:BVP:1980} and the references quoted therein. 

A different way to avoid the definition of a truncated boundary is to apply coordinate transforms.
Coordinate transforms have been applied successfully to the numerical solution of ordinary and partial differential equations on unbounded domains, see Grosch and Orszag \cite{Grosch:NSP:1977}, Koleva \cite{Koleva:NSH:2006} or Fazio and Jannelli \cite{Fazio:2014:FDS}.

This paper is concerned with two examples on the application of the \FBF \ of BVPs on a semi-infinite interval.
In both cases we are able to provide the exact solution of both the BVP and its \FBF .
Therefore, these problems can be used as benchmarks for the numerical methods applied to BVPs on a semi-infinite interval and to free BVPs. 
In this context, sometimes, it is possible to solve a given free BVP non-iteratively, see the survey by Fazio \cite{Fazio:1998:SAN}, whereas BVPs are usually solved iteratively.
Here, for two classes of free BVPs, we define non-iterative initial value methods which are referred in literature as non-iterative transformation methods (ITMs).
Indeed, non-ITMs can be defined within Lie's group invariance theory.
For the group invariance theory, the interested reader is referred to Bluman and Cole \cite{Bluman:1974:SMD}, Bluman and Kumei \cite{Bluman:1989:SDE}, Barenblatt \cite{Barenblatt:1996:SSI}, or Dresner \cite{Dresner:1999:ALT}.

Let us remark here that the first application of a non-iterative initial value method for the numerical solution of Blasius problem of boundary layer theory was given by T{\"o}pfer \cite{Topfer:1912:BAB}.
The algorithm devised by T{\"o}pfer for Blasius problem was redefined, and extended to a class of problems, within group invariance theory by Klamkin \cite{Klamkin:1962:TCB}.
The relationship between the invariance of the Blasius problem with respect to a linear group of transformations, the scaling group, and the applicability of a non-iterative initial value method was point out by Na \cite{Na:1967:TBC,Na:1968:FET}.
Moreover, in the same papers Na considered BVPs on finite intervals and the invariance with respect to a nonlinear group of transformations: the spiral group. 
Fazio and Evans, in \cite{Fazio:1990:SNA}, showed how to apply the scaling group to solve non-iteratively free BVPs.
The translation group of transformation was first used for the non-iterative solution of free BVPs by Fazio \cite{Fazio:1990:NVT}, see also the recent contribution by Fazio and Iacono \cite{Fazio:2010:TGN}.
Fazio \cite{Fazio:1992:MBH}, using the scaling group, defined a non-iterative initial value method for the numerical solution of a free boundary value problem governed by a system of first order differential equations.
For this author's knowledge no other Lie's group of point transformations has been applied to define a non-iterative initial value method.

In the past, the main drawback of non-ITMs was that they were considered not widely applicable: see the critical considerations by Fox, Erickson and Fan \cite{Fox:1969:LBL}, Meyer \cite[pp. 97-98]{Meyer:1973:IVM}, Na \cite[p. 137]{Na:1979:CME} or Sachdev \cite[p. 218]{Sachdev:1991:NOD}. 
In fact, the simplest way in order to verify if a non-ITM is applicable to a particular problem is to use an inspectional analysis as shown by Seshadri and Na \cite[pp. 157-168]{Seshadri:1985:GIE}, cf. also the discussion on inspectional analysis by Birkhoff \cite[pp. 99-103]{Birkhoff:1960:HSL}.

If we consider any possible extension of non-ITMs, then we have to recall the extension of scaling invariance involving physical parameters by Na \cite{Na:1970:IVM}, see also Na \cite[Chapters 8-9]{Na:1979:CME}.
Moreover, as shown by Fazio \cite{Fazio:1994:FSE,Fazio:1994:NTM,Fazio:1996:NAN}, within scaling invariance theory, it is possible to define an iterative extension of our approach that can be applied to the most general class of free BVPs.

\section{Free boundary formulation}
In this paper we provide two examples that support the idea that we can deal with BVPs  defined on a semi-infinite interval via their  \FBF .
In particular, in order to explain the main idea behind our \FBF , we consider the class of BVPs 
\begin{align}\label{eq:BVPs}
&{\ds \frac{d^2u}{dx^2}
+ f\left(x, u,\frac{du}{dx}\right)} = 0 \ , \qquad x \in [0, \infty) \nonumber \\[-1.5ex]
&\\[-1.5ex]
&{\ds u(0) = u_0 \ ,  \qquad u(\infty)}
= u_\infty \nonumber
\end{align}
where $f(\cdot,\cdot,\cdot)$ is a given function of its arguments, and $u_0$ and $u_\infty$ are given constants. 
If we can assume that the first derivative of $u(x)$ goes monotonically to zero at infinity, then we replace the problem (\ref{eq:BVPs}) with its \FBF 
\begin{align}\label{eq:BVPs:FBF}
&{\ds \frac{d^2u_\e}{dx^2}
+ f\left(x, u_\e,\frac{du_\e}{dx}\right)} = 0 \ , \qquad x \in [0, x_\e] \nonumber \\[-1.5ex]
&\\[-1.5ex]
&{\ds u_\e(0) = u_0 \ ,  \qquad u_\e(x_\e)}
= u_\infty \ , \quad \frac{du_\e}{dx} (x_\e) = \e \nonumber
\end{align}
where $x_\e$ is an unknown free boundary and $0 \le |\e| \ll 1$ is a parameter.

We have to remark here that monotonic properties of the solution, its first and second derivative have been demonstrated by Countyman and Kannan \cite{Countryman:1994:CNB}, for the class of problems in (\ref{eq:BVPs}) where $f$ depends exclusively on $u$.    

The following theorem provides, under suitable smoothness conditions, the order of convergence (and the uniform convergence) of the solution of
(\ref{eq:BVPs:FBF}) to the solution of (\ref{eq:BVPs}).

\medskip
\noindent
\begin{theorem}\label{Th:Conv} 
Suppose $ u_{\epsilon}(x) $ and $ {\frac{\partial u_{\epsilon}}{\partial \epsilon }} (x) $ are continuous functions with respect to $ \epsilon $ (and also with respect to $ x $ in the related free boundary domain $ [0, x_\epsilon ] $) and that $ | \epsilon_1 | < | \epsilon_2 | \Rightarrow [0, x_{\epsilon _2}] \subset
[0, x_{\epsilon _1}] $ at least in a non-empty interval including $ \epsilon = 0 $, then 
\begin{eqnarray*}
& || u_\epsilon (x) - u(x) || \leq K | \epsilon |
\end{eqnarray*}
where $ || \cdot || $ is the maximum norm on $[0, x_\e]$ and $ K $ is a positive constant independent on $ \epsilon $.
\end{theorem}

\medskip
\noindent
The proof of this Theorem can be obtained along the lines of the proof for the convergence Theorem stated in \cite{Fazio:1996:NAN} for a \FBF \ for a class of problems governed by a third order differential equation.

The \FBF \ allows us to embed a BVP in (\ref{eq:BVPs}) into a class of problems involving the control parameter $\e$.
When we solve the \FBF \ (\ref{eq:BVPs:FBF}) numerically, we can fix a very small value of $|\e|$ and apply a grid refinement to verify whether the numerical results agree within a prefixed number of significant digits. 
Also, it is possible to fix a step size and let $\e$ goes to zero and verify whether $u_\e(x) \rightarrow u(x)$ together with $x_\e \rightarrow \infty$. 
Usually, it suffices to take $|\e| \in \left\{10^{-1}\right.$, $10^{-2}$, $10^{-3}$, $10^{-4}$,
 $10^{-5}$, $10^{-6}$, $\left. \dots \right\}$
and compare the obtained numerical results.
Let us remark here that sometimes it is possible to solve the \FBF \ non-iteratively, see the survey by Fazio \cite{Fazio:1998:SAN}, whereas the BVP (\ref {eq:BVPs}) is usually solved iteratively.
 
\section{Two examples for the \FBF}
\begin{figure}[!hbt]
\centering
\psfrag{x}[][]{$x$} 
\psfrag{u}[][]{$u(x)$} 
\psfrag{p01}[][]{$P = 0.1$} 
\psfrag{p1}[][]{$P = 1$} 
\psfrag{p10}[][]{$P = 10$} 
\includegraphics[width=\www]{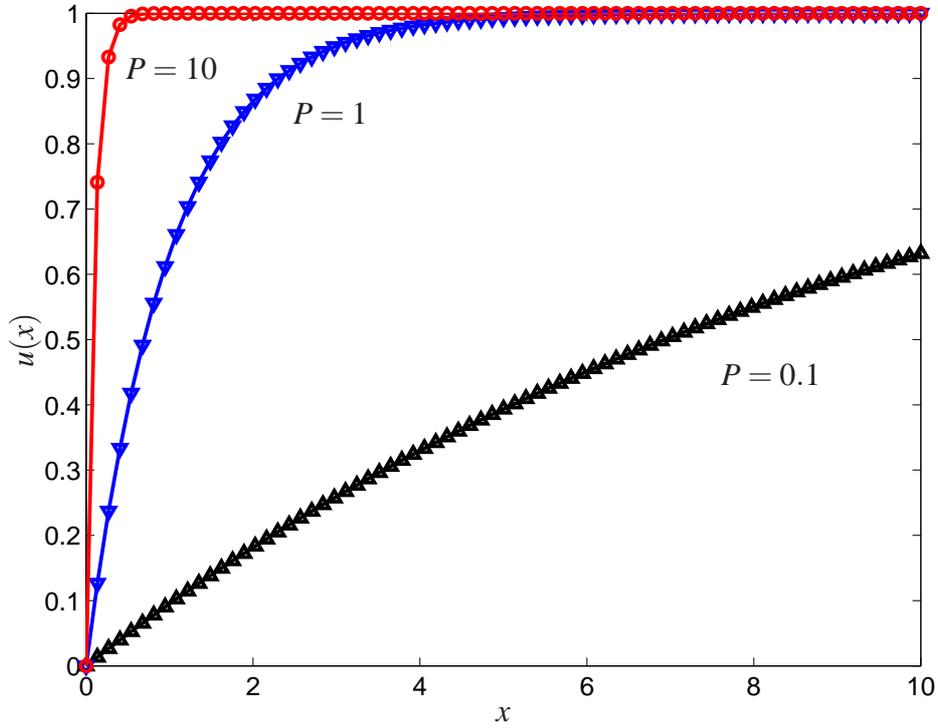}
\caption{The solution (\ref{eq:sol1}) for different values of  $P$. The symbols stand for:
$\circ$ $P=10$, $\triangledown$ $P=1$, and $\vartriangle$ $P=0.1$.}
\label{fig:BVPes1}
\end{figure}

As a first example we consider the linear problem
\begin{align}\label{eq:ex1}
& {\ds \frac{d^2u}{dx^2} + P \frac{du}{dx}} = 0 \ , \qquad x \in [0, \infty) \nonumber \\[-1.5ex]
& \\[-1.5ex]
& u(0) = 0 \ , \qquad u(\infty) = 1 \nonumber 
\end{align}
where $P$ is a positive constant.
The solution of (\ref{eq:ex1}) is easily found to be
\begin{equation}\label{eq:sol1}
u(x) = 1-e^{-Px}
\end{equation}
so that the missing initial condition is
equal to $P$, that is $\frac{du}{dx}(0) = P$.
Figure \ref{fig:BVPes1} shows the solution (\ref{eq:sol1}) of the BVP (\ref{eq:ex1}) for different values of  $P$.
The bigger is the value of $P$, the harder is to solve the BVP numerically.
In fact, for large values of $P$ the solution has a boundary layer near $x=0$. 

\begin{figure}[!hbt]
\centering
\psfrag{x}[][]{$x$} 
\psfrag{u}[][]{$u(x), u_\e(x)$} 
\psfrag{P}[][]{$P = 1$} 
\includegraphics[width=\www]{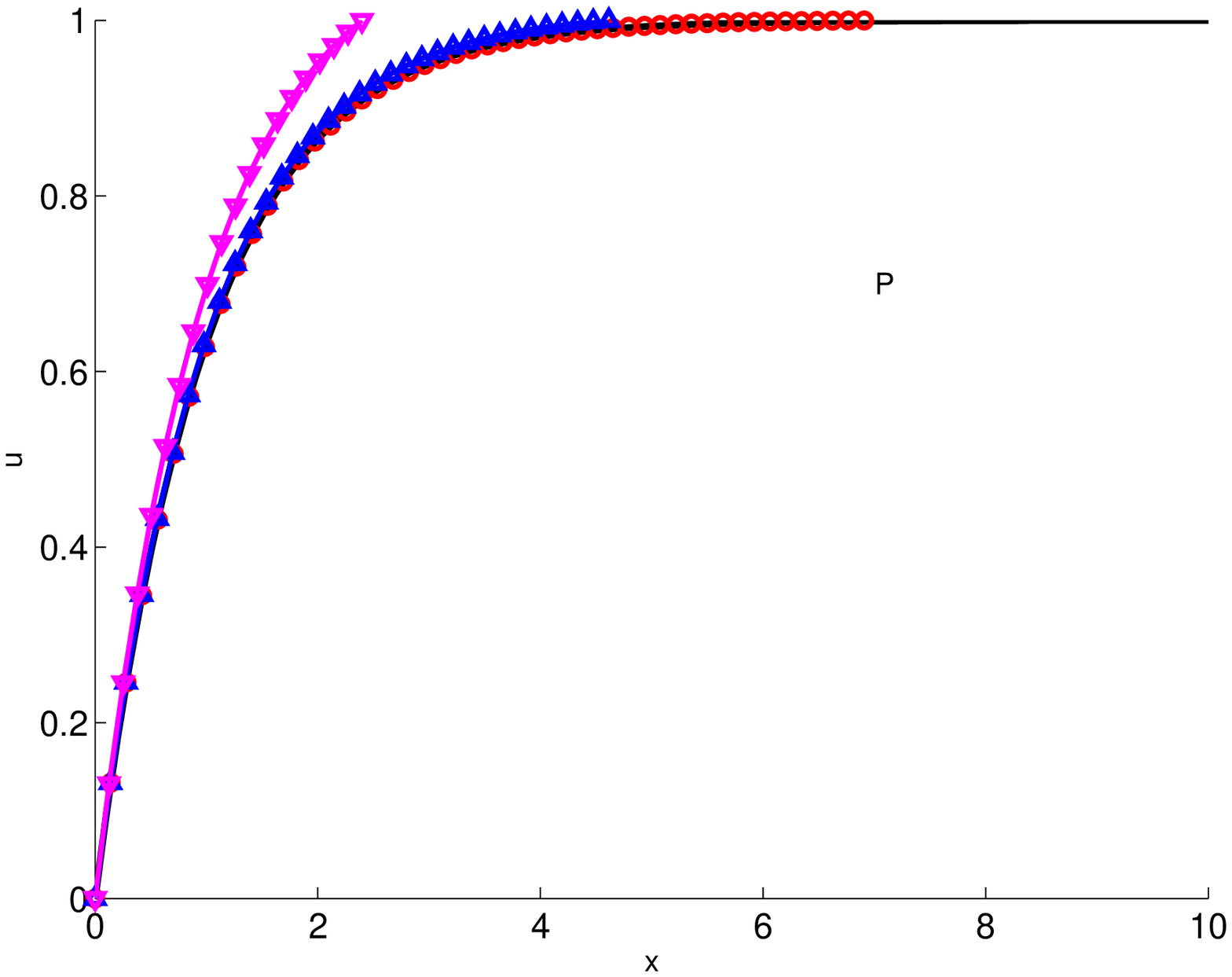}
\caption{The solution (\ref{eq:solfbf1}) for $P = 1$ and different values of  $\e$. The symbols stand for:
$\mathbf{-}$ the exact solution, $\triangledown$, $\vartriangle$  and $\circ$ the free boundary solution $u_\e$ with $\e=0.1$, $\e=0.01$ and $\e=0.001$, respectively.}
\label{fig:FBFes1}
\end{figure}

Let us consider now the \FBF \ for (\ref{eq:ex1})
\begin{align}\label{eq:ex1:fbf}
& {\ds \frac{d^2u_\e}{dx^2} + P \frac{du_\e}{dx}} = 0 \ , \qquad x \in [0, x_\e] \nonumber \\[-1.5ex]
& \\[-1.5ex]
& u_\e(0) = 0 \ , \qquad u_\e(x_\e) = 1 \ , \qquad \frac{du_\e}{dx}(x_\e) = \e \ , \nonumber 
\end{align}
with $0 \le \e \ll 1$.
The solution of (\ref{eq:ex1:fbf}) is given by
\begin{equation}\label{eq:solfbf1}
u_\e (x) = {\ds \frac{P+\e}{P}\left(1-e^{-Px}\right)} \ , \quad x_\e = - {\ds \frac{1}{P} \ln\left(\frac{\e}{P+\e}\right)} \ .
\end{equation}
Therefore, we can easily verify that as $\e$ goes to zero the solution $u_\e(x)$ of the \FBF \
(\ref{eq:ex1:fbf}) converges to the solution $u(x)$ of the original problem (\ref{eq:ex1}) and the free boundary $x_\e$ goes to infinity.
Moreover, we realize that the obtained approximation becomes the more accurate the more $\e$ is near zero, see figure \ref{fig:FBFes1}.

Let us remark here that the same exact solutions (\ref{eq:sol1}) and (\ref{eq:solfbf1}) are still valid if we replace the governing differential equation, in the BVP (\ref{eq:ex1}) and its \FBF \ (\ref{eq:ex1:fbf})
with the non-autonomous one
\begin{equation}
{\ds \frac{d^2u}{dx^2} + P^2 e^{-P x}} = 0 \ ,
\end{equation}
where we substitute $u = u_\e$ in the free boundary case.

Replacing a linear problem with a nonlinear one can be justified, from a numerical viewpoint, only by considering that in this way we overcome the singularity related to the boundary condition prescribed at infinity.
Of course, when the original problem is a nonlinear one a \FBF \ for it can be really convenient to solve numerically.

\begin{figure}[!hbt]
\centering
\psfrag{x}[][]{$x$} 
\psfrag{u}[][]{$u(x)$} 
\psfrag{p01}[][]{$P = 0.1$} 
\psfrag{p1}[][]{$P = 1$} 
\psfrag{p10}[][]{$P = 10$} 
\includegraphics[width=\www]{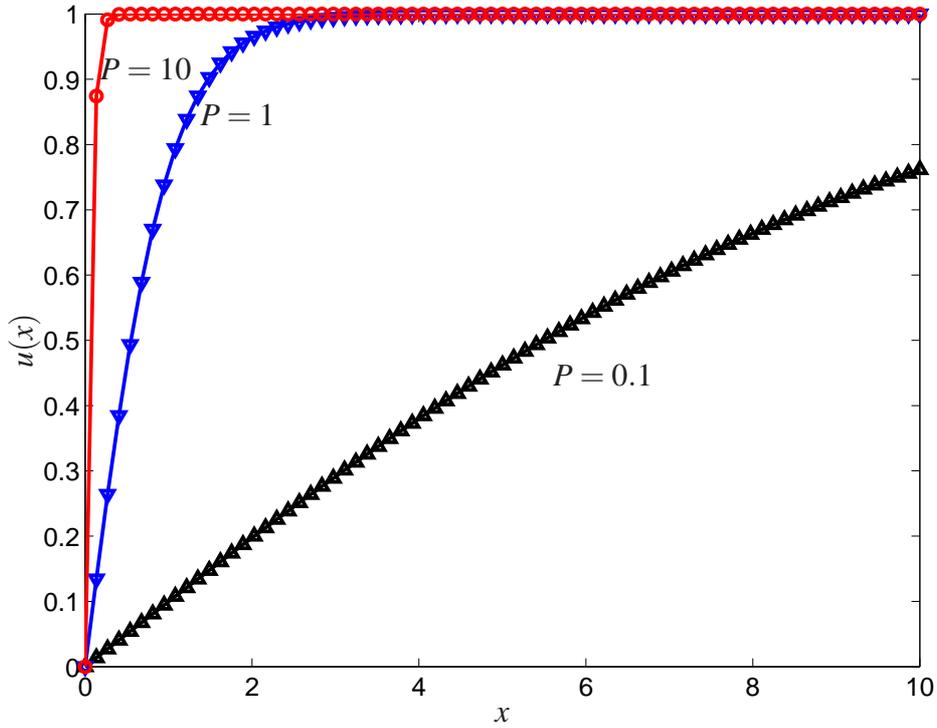}
\caption{The solution (\ref{eq:sol2}) for different values of  $P$. The symbols stand for:
$\circ$ $P=10$, $\triangledown$ $P=1$, and $\vartriangle$ $P=0.1$.}
\label{fig:BVPes2}
\end{figure}

As a second example we consider the nonlinear problem
\begin{align}\label{eq:ex2}
& {\ds \frac{d^2u}{dx^2} + 2 P u \frac{du}{dx}} = 0 \ , \qquad x \in [0, \infty) \nonumber \\[-1.5ex]
& \\[-1.5ex]
& u(0) = 0 \ , \qquad u(\infty) = 1 \ , \nonumber 
\end{align}
where, again, $P$ is a positive constant.
The solution of (\ref{eq:ex2}) is given by
\begin{equation}\label{eq:sol2}
u(x) = \tanh(Px) \ ,
\end{equation}
and, again, $\frac{du}{dx}(0) = P$.
Figure \ref{fig:BVPes2} shows the solution (\ref{eq:sol2}) of the BVP (\ref{eq:ex2}) for different values of  $P$.
Again, for large values of $P$ the solution has a boundary layer near $x=0$. 
It can be easily verified that, for instance by comparing figure (\ref{fig:BVPes1}) with figure (\ref{fig:BVPes2}), for the same value of the parameter $P$, the BVP (\ref{eq:ex2}) is more challenging than the BVP (\ref{eq:ex1}). 

\begin{figure}[!hbt]
\centering
\psfrag{x}[][]{$x$} 
\psfrag{u}[][]{$u(x), u_\e(x)$} 
\psfrag{P}[][]{$P = 1$} 
\includegraphics[width=\www]{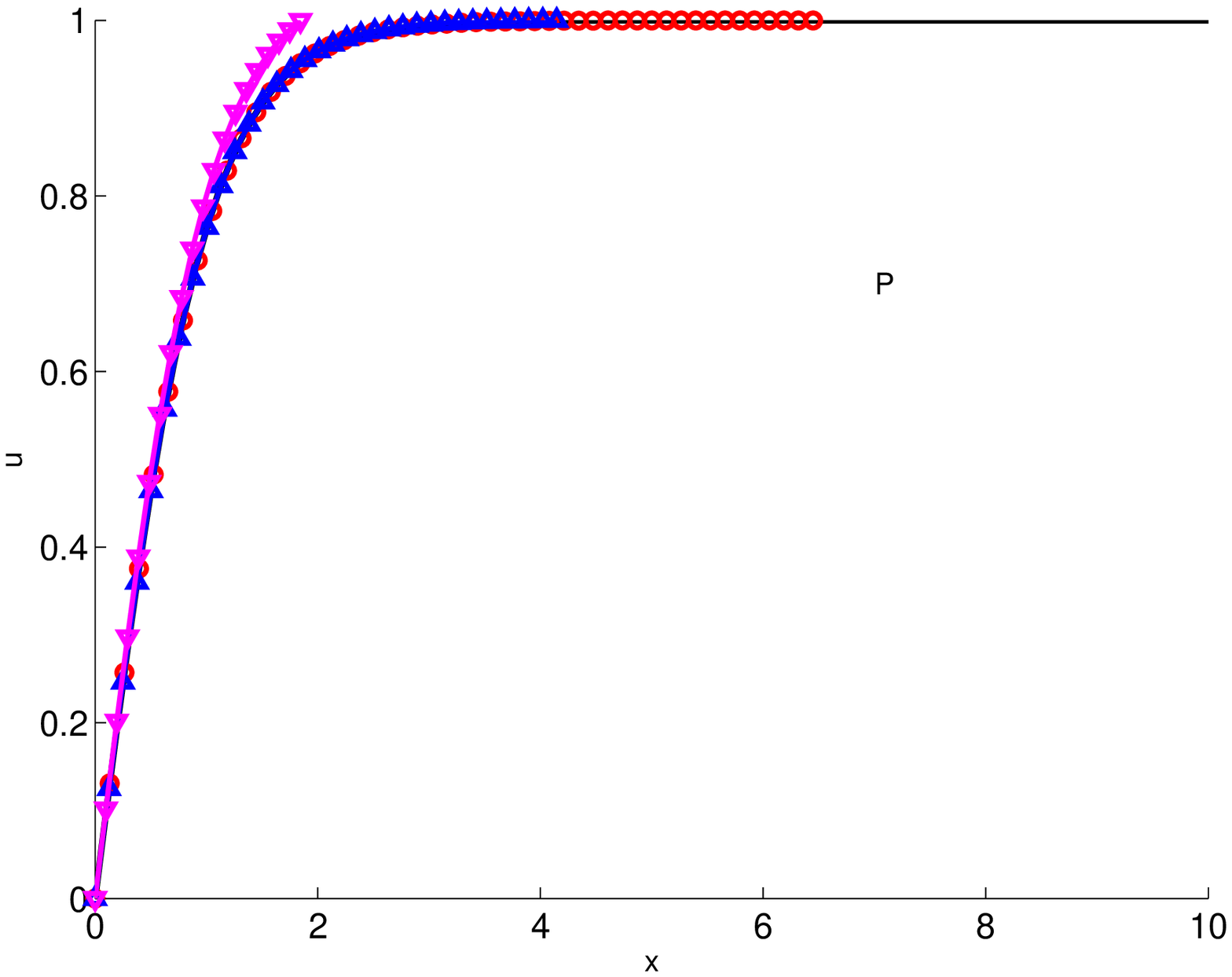}
\caption{The solution (\ref{eq:solfbf2}) for $P = 1$ and different values of  $\e$. The symbols stand for:
$\mathbf{-}$ the exact solution, $\triangledown$, $\vartriangle$  and $\circ$ the free boundary solution $u_\e$ with $\e=0.1$, $\e=0.001$ and $\e=0.00001$, respectively.}
\label{fig:FBFes2}
\end{figure}

Let us consider now the \FBF \ for (\ref{eq:ex2})
\begin{align}\label{eq:ex2:fbf}
& {\ds \frac{d^2u_\e}{dx^2} + 2 P u_\e \frac{du_\e}{dx}} = 0 \ , \qquad x \in [0, x_\e] \nonumber \\[-1.5ex]
& \\[-1.5ex]
& u_\e(0) = 0 \ , \qquad u_\e(x_\e) = 1 \ , \qquad {\ds \frac{du_\e}{dx}(x_\e)} = \e \ , \nonumber 
\end{align}
with $0 \le \e \ll 1$.
The positive solution of (\ref{eq:ex2:fbf}) is given by
\begin{equation}\label{eq:solfbf2}
u_\e (x) = {\ds - \frac{1}{C}\tanh{(Px)}} \ , \quad x_\e = {\ds \frac{1}{2P} \ln\left(\frac{1-C}{1+C}\right)} \ ,
\end{equation}
where ${\ds C = \left(\e - \sqrt{\e^2+4 P^2}\right)/{2P}}$. 
Also in this case, as $\e$ goes to zero the solution $u_\e(x)$ of the \FBF \
(\ref{eq:ex2:fbf}) converges to the solution $u(x)$ of the original problem (\ref{eq:ex2}) and the free boundary $x_\e$ goes to infinity.
Moreover, also in this case the obtained approximation becomes the more accurate the more $\e$ is close to zero,
figure \ref{fig:FBFes2}.

\section{Non-ITMs for the \FBF}
As we mentioned in the introduction numerical TMs can be defined within Lie's group invariance theory.

\subsection{Spiral and translation invariance}
Let us define a non-ITM for the class of free BVPs
\begin{align}\label{eq:TraClass}
&\frac{d^2u}{dx^2}= u \; \Omega\left(u e^{-\omega x},\frac{du}{dx} e^{-\omega x}\right) \ , \qquad x \in [0, s] \nonumber\\
&u(0)  = \alpha \ , \\
&u(s)= \beta \; e^{\omega s}\quad , \quad \frac{du}{dx}(s)=\gamma \; e^{\omega s} \ ,\nonumber
\end{align}
where $ \Omega (\cdot, \cdot) $ is an arbitrary function of its arguments, $\alpha$, $\beta$, $\gamma$ and $\omega$ are arbitrary constants, and $s>0$ is the unknown free boundary.
The governing differential equation and the two free boundary conditions are left invariant by the spiral group
\begin{equation}\label{t1}
x^*=x+\mu\quad , \quad s^*=s+\mu\quad , \quad u^*=e^{\omega 	\mu}u \ ,
\end{equation}
where $\mu$ is the group parameter.

Let us remark here that by setting $\omega = 0$ we can recover from (\ref{eq:TraClass}) the class of problems invariant with respect to the translation group of transformation defined by (\ref{t1}) for $\omega = 0$. 

Using the transformation (\ref{t1}) we can define the following non-iterative algorithm for the numerical solution of (\ref{eq:TraClass}):
\begin{enumerate}
\item Input $s^*$. 
\item Solve backwards from $ s^* $ to $ x_{\alpha}^* $ (to be defined below) the following auxiliary IVP 
\begin{align}\label{eq:IVP1}
&\frac{d^2u^*}{dx^{*2}}= u^* \; \Omega\left(u^* e^{-\omega x^*},\frac{du^*}{dx^*} e^{-\omega x^*}\right) \ , \qquad x^* \in [x_\alpha^*, s^*] \nonumber \\[-1ex]
& \\[-1ex]
&u^*(s^*)=\beta \; e^{\omega s^*}\ , \qquad \frac{du^*}{dx^*}(s^*)=\gamma \; e^{\omega s^*}\ ,	\nonumber
\end{align}
using an {\it event locator} in order to find $ x_{\alpha}^* $ such that 
\begin{equation}\label{eq:evloc} 
u^*(x_{\alpha}^*)  = \alpha \ .
\end{equation} 
\item Compute the group parameter
\begin{equation}\label{eq:mu}
\mu = x_{\alpha}^* \ ,
\end{equation} 
the free boundary
\begin{equation}\label{eq:freeb}
s = s^*-\mu \ ,
\end{equation} 
and the missing initial condition
\begin{equation}\label{eq:mics}
\frac{du}{dx}(0) = e^{-\omega \mu} \frac{du^*}{dx^*}(x_{\alpha}^*) \ .
\end{equation}
\item
Compute the transformed solution
\begin{equation}\label{eq:trasol}
u(x) = e^{-\omega \mu} u^*(x^*) \ , \qquad \frac{du}{dx}(x) = e^{-\omega \mu} \frac{du^*}{dx^*}(x^*) \ .
\end{equation}
where, of course, $x = x^*-\mu$. 
\end{enumerate}

We define now a simple event locator which is suitable to be applied with the non-ITM for the numerical solution of (\ref{eq:TraClass}).
Let us consider the case $\alpha < \beta e^{\omega s^*}$, the case $\alpha > \beta e^{\omega s^*}$ can be treated in a similar way.
We integrate the auxiliary IVP (\ref{eq:IVP1}) until we get at a mesh point $x_k^*$ where $u^*(x_k^*) < \alpha$, and we compute  
\begin{equation}\label{eq:evloc1}
x_\alpha^* = x_k^* + (\alpha - u_k^*) \frac{\Delta x^*}{u_k^*-u_{k-1}^*} \ .
\end{equation}
Then, we repeat the last step with the smaller step size given by
\begin{equation}\label{eq:evloc2}
\Delta x_\alpha^* = x_\alpha^* - x_k^* \ .
\end{equation}
In defining the last step size in equation (\ref{eq:evloc1})-(\ref{eq:evloc2}), we apply a first order (linear) Taylor formula
at $x_k^*$ where we have replaced the first derivative with a backward finite difference approximation. 

Let us notice here that this non-ITM generalize the one proposed by Fazio and Iacono \cite{Fazio:2010:TGN} for the numerical solution of free boundary problems with governing equations invariant with respect to a translation group.

Moreover, under suitable hypotheses we can define, within group invariance theory, a transformation of variables allowing us to rewrite each free BVPs belonging to (\ref{eq:TraClass}) with $\omega = 0$ as a problem in the class of free BVPs (\ref{eq:ScaClass}) that will be considered in the next subsection, see \cite{Fazio:1990:NVT} for the details.

\subsection{Scaling invariance}
Let us define a non-ITM for the class of free BVPs
\begin{align}\label{eq:ScaClass}
& {\displaystyle \frac{d^{2}u}{dx^2}} = u^{1-2\delta}
\Phi \left( x u^{-\delta},
{\displaystyle \frac{du}{dx}} u^{\delta-1}
\right) \ , \qquad x \in [0, s] \nonumber \\
& u(0) = 0 \ , \\
& u(s) =  \beta \ ,
\qquad
{\displaystyle \frac{du}{dx}}(s) = \gamma \ , \nonumber
\end{align}
where $ \Phi (\cdot, \cdot) $ is an arbitrary function of its arguments,
$ \delta \neq 0$, $\beta \neq 0$ and $\gamma$ are arbitrary constants, and $s>0$ is the unknown free boundary.
The governing differential equation and the boundary condition at $x=0$ are invariant with respect to the scaling group:
\begin{equation}\label{eq:sgroup}  
x^* = \lambda ^\delta x, \qquad s^* = \lambda ^\delta s,
\qquad u^* = \lambda u  \ ,
\end{equation}
where $\lambda$ is the group parameter.
Using the invariance properties, we can define the following non-iterative algorithm for the numerical solution of (\ref{eq:ScaClass}):
\begin{enumerate}
\item Input $s^* > 0$, $v_0 \gg 1$ and $0 < \tau \ll 1$. 
\item Solve the auxiliary IVP
\begin{align}\label{eq:ex2:ivp}
& {\displaystyle \frac{d^{2}u^*}{dx^{*2}}} = u^{*1-2\delta}
\Phi \left( x u^{*-\delta}, {\displaystyle \frac{du^*}{dx^*}} u^{*\delta-1}\right) \ , \qquad x \in [0, s^*]  \nonumber \\[-1.5ex]
& \\[-1.5ex]
& u^*(0) = 0 \ , \qquad {\ds \frac{du^*}{dx^*}(0)} =  v_0 \ . \nonumber 
\end{align}
\item Compute $\lambda$ by 
\begin{equation}\label{eq:lambda1}
\lambda = \frac{u^*(s^*)}{\beta} \ ,
\end{equation}
and the free boundary
\begin{equation}\label{eq:FB2}
s = \lambda^{-\delta} s^*
\end{equation}
\item Rescale the numerical solution to get $u(x)$ and $\frac{du}{dx}(x)$ according to (\ref{eq:sgroup}).
In particular, we find
\begin{align}
& {\ds \frac{du}{dx}(s)} = \lambda^{\delta-1} {\ds \frac{du^*}{dx^*}(s^*)} \nonumber \\[-1ex]
& \\[-1ex]
& {\ds \frac{du}{dx}(0)} = \lambda^{\delta-1} v_0 \ . \nonumber 
\end{align}
\end{enumerate}
If the computed value of $\frac{du}{dx}(s)$ is bigger than $\tau$, then set a lager value for $\ni_0$ and repeat the computation. 
It is evident that there is no need to rewrite the given free BVP in standard form, as suggested by Ascher and Russell \cite{Ascher:1981:RBV}, because we can choose $s^*$ at our convenience and therefore we usually set $s^*=1$.   

The algorithm presented above is an original variant of the non-ITM defined in \cite{Fazio:1990:SNA}, where the boundary conditions in (\ref{eq:ScaClass}) were replaced by
\begin{equation}\label{eq:BCs}
u(0) = \alpha \ , \qquad u(s) =  \beta s^{1/\delta}\ , \quad {\displaystyle \frac{du}{dx}}(s) = \gamma s^{(1-\delta)/\delta}\ ,
\end{equation}
where $ \alpha \neq 0$, $\beta$ and $\gamma$ are arbitrary constants and we had to integrate backwards on $[0, s^*]$ because these free boundary values are invariant, but the boundary condition at zero is not invariant, under the scaling group (\ref{eq:sgroup}).

\section{Numerical tests}
In this section we report on the application of the two non-ITMs defined in the previous section.
To this end we solve numerically the \FBF s defined in section 3.

\subsection{Translation invariance}
As an example we consider here the non-iterative numerical solution of the \FBF \ (\ref{eq:ex1:fbf}) with $P = 1$, $ \e = 1\mbox{D}-06$.
With a simple change of variables, the problem (\ref{eq:ex1:fbf}) clearly belongs to (\ref{eq:TraClass}) with  $\omega = 0$.
In table \ref{tab:RKes1} we list sample numerical results. 
We used a uniform grid, with $\Delta x^* = -0.05$ and $x_\e^* =1$, and applied 
the classical fourth order Runge-Kutta method (RK4) \cite{Runge:1895:UNA,Kutta:1901:BNI}, the sixth order Runge-Kutta method (RK6), and the eighth order Runge-Kutta method (RK8) as reported by Butcher in \cite{Butcher:NMO:2003} on page 178 and 180, respectively.
 \input{TABes1}
Since the numerical results obtained by RK6 and RK8 are very close, we can infer that in order to improve the numerical accuracy we need to consider a grid refinement. 

\begin{figure}[!hbt]
\centering
\psfrag{x}[][]{$x$} 
\psfrag{u}[][]{$u_\e(x)$, $\frac{du_\e}{dx}(x)$} 
\includegraphics[width=.8\textwidth]{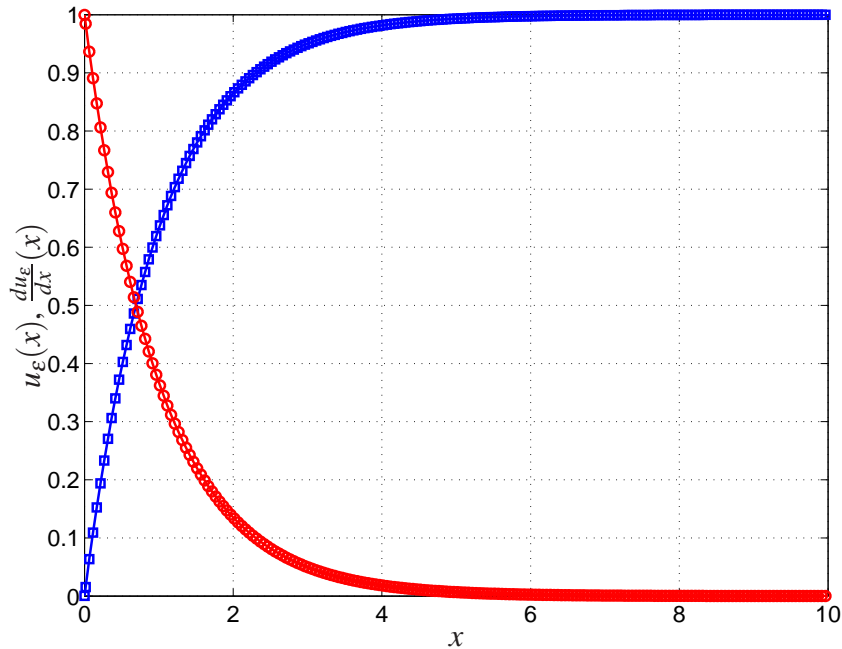}
\caption{Numerical solution of the \FBF \ (\ref{eq:ex1:fbf}). The symbols stand for: $\square$ $\frac{du_\e}{dx}(x)$, and $\circ$ $u_\e(x)$. We notice that $x_\e \approx 13.8$.}
\label{fig:Es1}
\end{figure}

In table \ref{tab:RKes1:Dx} we list the numerical results obtained by RK6 and a grid refinement with the reported step sizes. 
 \input{TABes1Dx}
The reported numerical results clearly indicate that we are able to get an accurate numerical approximation of the BVP (\ref{eq:ex1}) non-iteratively.
This can be also realized by comparing the exact solution, for $P=1$, plotted on figure \ref{fig:BVPes1} 
with the numerical solution shown on figure \ref{fig:Es1}.
Figure \ref{fig:Es1} is a frame, on $x\in [0, 10]$, of the numerical solution, computed by the RK6 solver with $\Delta x^* = -0.05$.
From figure \ref{fig:Es1} we can notice how the last step, which is the one close to the origin, is smaller with respect to the previous ones.
This is due to the application of our simple event locator defined by equations (\ref{eq:evloc1})-(\ref{eq:evloc2}).

In table \ref{tab:RKes1:eps} we report the numerical results obtained by fixing a value for the step size and repeating the computation for several decreasing values of $\e$.
\input{TABes1eps}
As it is easily seen, we can verify numerically that, the smaller the value of $\e$ the larger the free boundary value $x_\e$. 

\subsection{Scaling invariance}
As an example we consider here the non-iterative numerical solution of the \FBF \ (\ref{eq:ex2:fbf}). 
With a simple change of variables, the problem (\ref{eq:ex2:fbf}) clearly belongs to (\ref{eq:ScaClass}).
The governing differential equation is invariant on condition that $ 1-2\delta = 2-\delta$, that is $\delta = -1$.

\begin{figure}[!hbp]
\centering
\psfrag{x}[][]{$x$} 
\psfrag{u}[][]{$u_\e(x)$, $\frac{du_\e}{dx}(x)$} 
\psfrag{x*}[][]{$x^*$} 
\psfrag{u*}[][]{$u_\e^*(x^*)$, $\frac{du_\e^*}{dx^*}(x^*)$} 
\includegraphics[width=.8\textwidth]{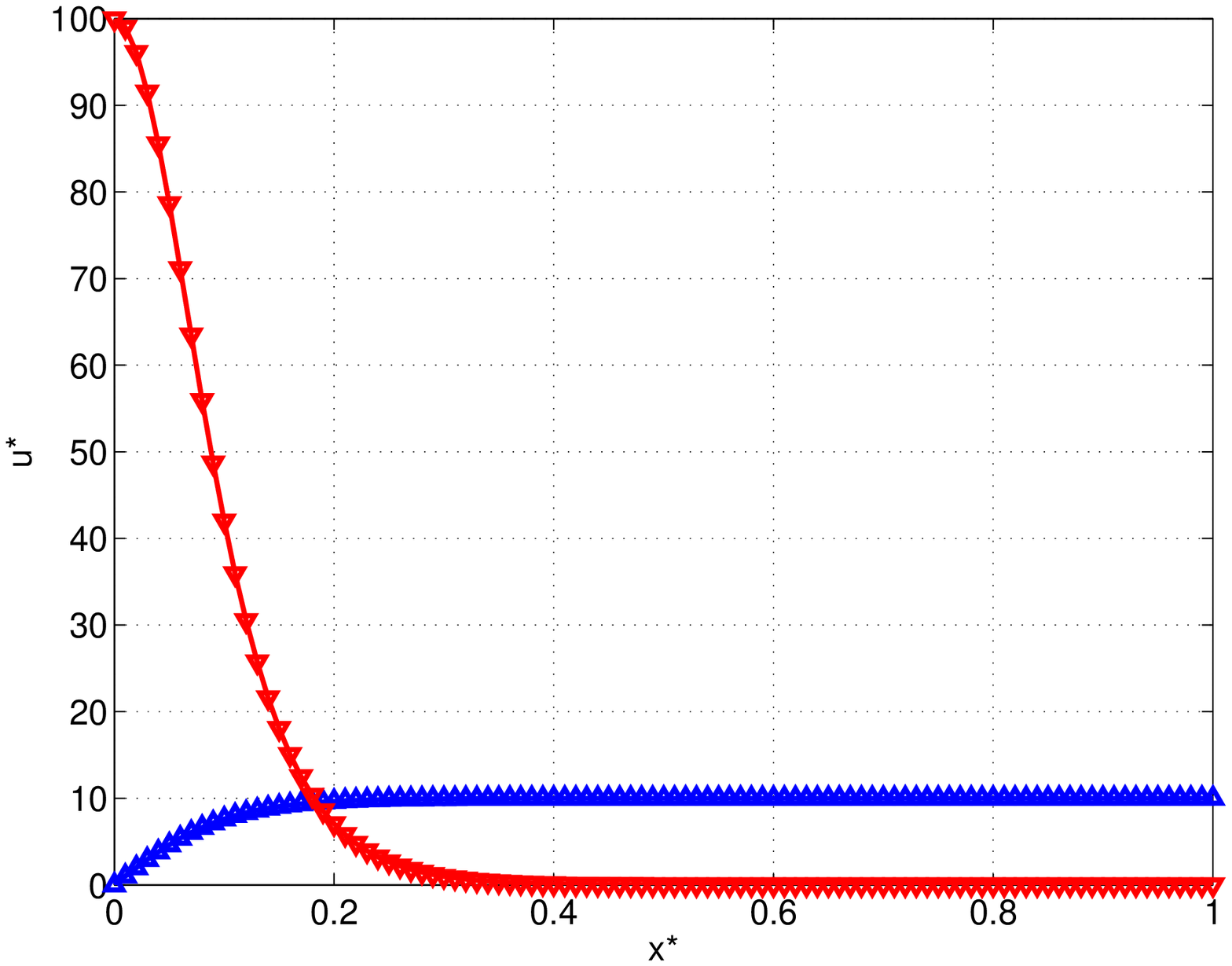} \\
\includegraphics[width=.8\textwidth]{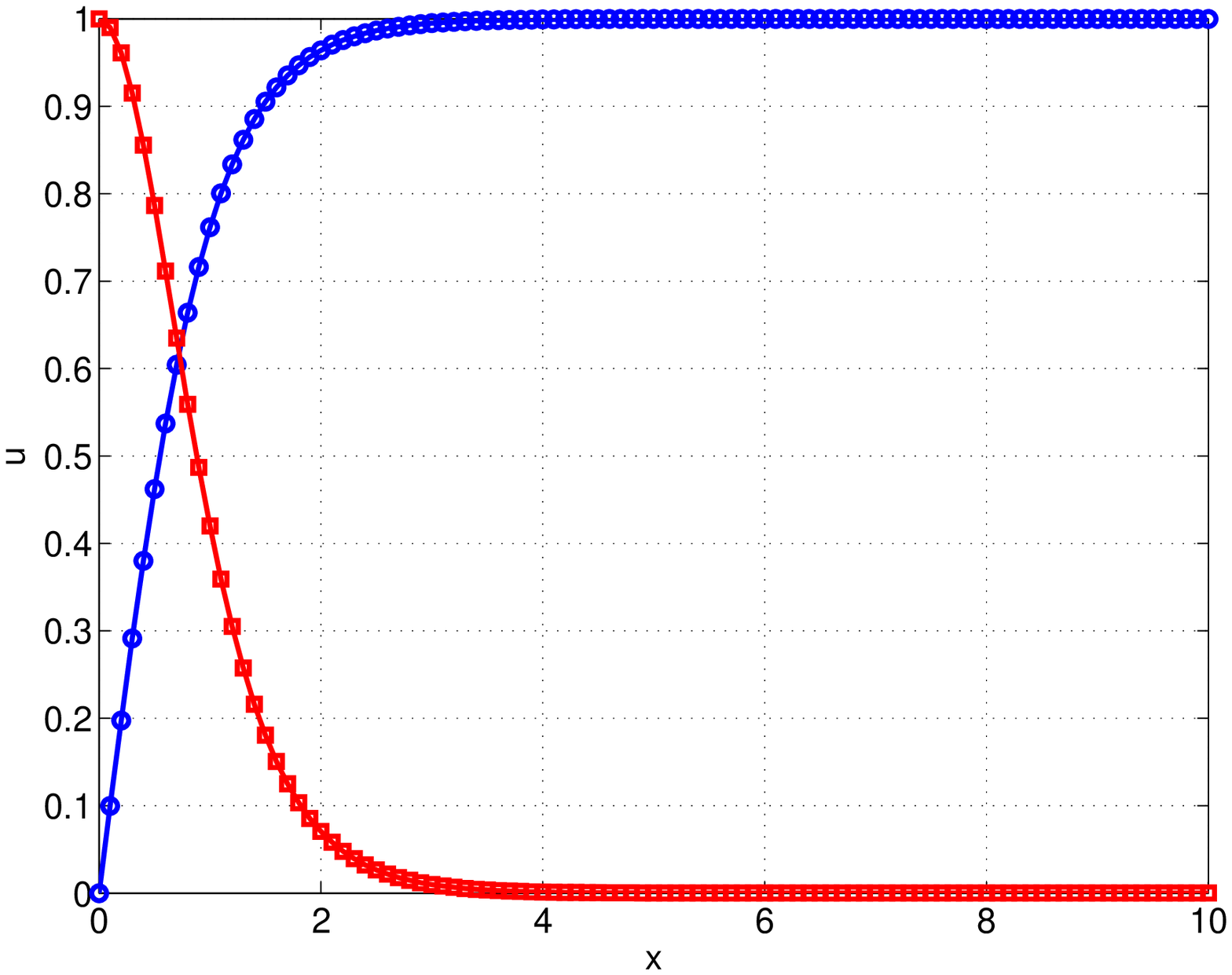}
\caption{Top frame: solution of the initial value problem (\ref{eq:ex2:ivp}). Bottom frame solution of the \FBF \ 
(\ref{eq:ex2:fbf})  with $v_0 = 100$. The symbols stand for:
$\triangledown$ $\frac{du_\e^*}{dx^*}(x^*)$, $\vartriangle$ $u_\e^*(x^*)$, $\square$ $\frac{du_\e}{dx}(x)$, and $\circ$ $u_\e(x)$.}
\label{fig:Es2}
\end{figure}

In table \ref{tab:RKes2} we list the numerical results obtained by setting $P = 1$, $x_\e^* =1$, $v_0 = 100$ and $\tau = 1 \; \mbox{D}-06$.
We used a uniform grid with $\Delta x^* = 0.01$ and the same Runge-Kutta solvers of the previous subsection.
 \input{TABes2}
The reported numerical results clearly indicate that we are able to get an accurate numerical approximation of the BVP (\ref{eq:ex2}) non-iteratively.
This can be also realized by comparing the exact solution, for $P=1$, plotted on figure \ref{fig:BVPes2} 
with the numerical solution shown on the bottom frame of  figure \ref{fig:Es2}.
In figure \ref{fig:Es2} we plot the numerical solution of the initial value problem (\ref{eq:ex2:ivp}), obtained by RK6, as well as the rescaled solution for the \FBF \ of the BVP (\ref{eq:ex2}).

\begin{figure}[!hbp]
\centering
\psfrag{x*}[][]{$x^*$} 
\psfrag{u*}[][]{$u_\e^*(x^*)$, $\frac{du_\e^*}{dx^*}(x^*)$} 
\includegraphics[width=.8\textwidth]{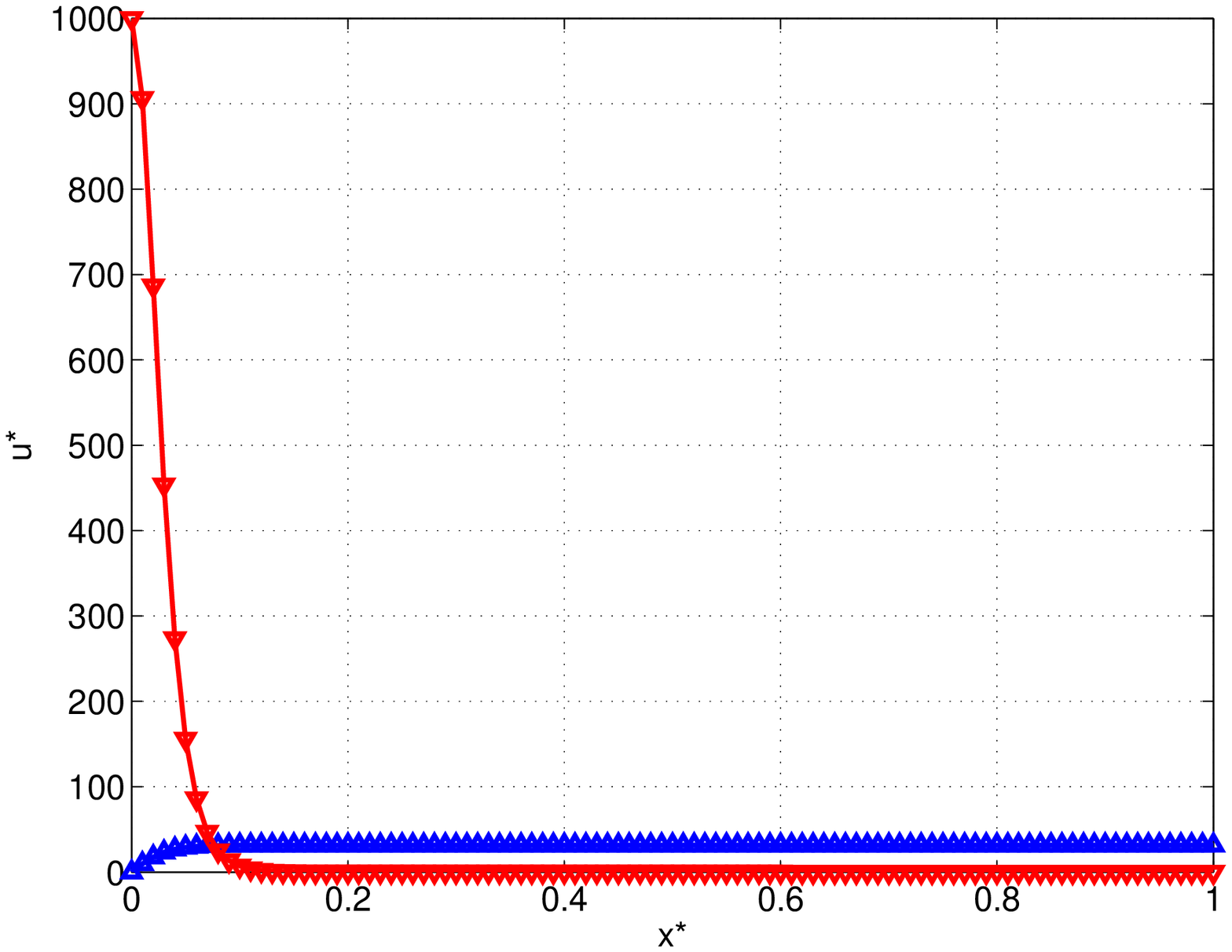} 
\caption{Numerical solution of the initial value problem (\ref{eq:ex2:ivp}) with $v_0 = 1000$. The symbols stand for:
$\triangledown$ $\frac{du_\e^*}{dx^*}(x^*)$ and $\vartriangle$ $u_\e^*(x^*)$.}
\label{fig:Es2:1000}
\end{figure}

We would like to remark that, in order to solve the \FBF \ non-iteratively we can use a whatever large value of $v_0$.
For instance, by setting $v_0 = 1000$ and again $\Delta x^* = 0.01$ we have computed with RK6 the values ${ \frac{du_\e}{dx}(x_\e)} \approx 1.37\mbox{D}-27$, $x_\e \approx 31.62$, and ${\frac{du_\e}{dx}(0)} \approx 0.999978$.
Indeed, we got a less accurate value for the first derivative of $u_\e(x)$ at $x=0$, and this depends on the use of the same grid on $[0,1]$ although we have to cope with a faster transitory of $\frac{du_\e^*}{dx^*}(x^*)$ near the origin, see figure \ref{fig:Es2:1000}.
However, we verify numerically that the smaller the value of $\e$ the larger the free boundary value $x_\e$. 

\section{Concluding remarks}
In this paper we have reported on the \FBF \ for BVPs defined on a semi-infinite interval. 
In this context, we discussed two simple test problems where we are able to get the exact analytical solution of both the original BVPs and of their \FBF . 
Therefore, these problems can be used as benchmarks for numerical methods applied to BVPs on a semi-infinite interval and to free BVPs. 
Moreover, we emphasized how for two classes of free BVPs we can define non-iterative initial value methods, whereas BVPs are usually solved iteratively.
These non-iterative methods can be deduced within Lie's group invariance theory.
We also applied the non-iterative methods to the two introduced \FBF s and we obtained meaningful numerical results.
 
Let us indicate a few problems in the applied sciences where our non-ITMs can be used.
For instance, the non-ITM defined using the translation group can be applied to the \FBF \ of the BVP 
\begin{align}\label{eq:bvp:Colloids}
&\frac{d^2u}{dx^2}-2 \; \sinh(u) = 0  \ , \qquad x \in [0, \infty) \nonumber \\[-1ex]
&\\[-1ex]
& u(0) = c \ , \qquad u(\infty)=0 \ , \nonumber 
\end{align}
where $c$ is a positive constant, arising in the theory of colloids, see Countryman and Kannan \cite{Countryman:1994:CNB} and the references quoted therein.
Moreover, the non-ITM defined using the scaling group can be applied to the \FBF \ of the BVP, studied by Seshadri and Na \cite{Seshadri:1975:ISN}, 
\begin{align}\label{eq:bvp:SNa}
&\frac{d^2u}{dx^2}+ m x \left(\frac{du}{dx}\right)^{(2-q)} = 0 \ , \qquad x \in [0, \infty) \nonumber \\[-1ex]
&\\[-1ex]
& u(0) = 0 \ , \qquad u(\infty)=1 \ , \nonumber 
\end{align}
where $m$ and $q$ are physical parameters, arising in the study of longitudinal impact to a nonlinear viscoplastic thin rod.
We would like to mention that the non-ITMs defined above can be easily generalized to the \FBF \ of BVPs governed by systems of differential equations like those derived by Dresner \cite[pp. 69-97]{Dresner:1983:SSN} within the study of wave propagation problems.

Finally, let us report on a possible way to extend the non-iterative methods defined in this paper to free boundary problems that are invariant with respect to a generic Lie group.
To this end we assume that a given problem is invariant under the one-parametr group of point transformations
\begin{equation}\label{eq:Lieg}
x^* = x + \lambda X(x, u) \ , \qquad u^* = u+\lambda U(x, u) \ ,
\end{equation}
where $\lambda$ is the group parameter and $X(x, u)$ and $U(x, u)$ are the so-called group generators.
The infinitesimal generator $D$ of (\ref{eq:Lieg}) is given by
\begin{equation}\label{eq:inf1}
D = X(x, u) \frac{\partial}{\partial x} + U(x, u) \frac{\partial}{\partial u} \ .
\end{equation}
If we introduce the variable transformation $y = y(x, u)$ and $w =w(x, u)$, then we can rewrite the infinitesimal generator (\ref{eq:inf1}) as
\begin{equation}\label{eq:inf2}
D = \left(X \frac{\partial y}{\partial x} + U \frac{\partial y}{\partial u}\right) \frac{\partial}{\partial y} + \left(X \frac{\partial w}{\partial x} + U \frac{\partial w}{\partial u}\right) \frac{\partial}{\partial w} \ .
\end{equation}
Now we can choose between two possible alternative: the so-called canonical variables, see Bluman and Kumei \cite{Bluman:1989:SDE}, that transforms as
\begin{equation}\label{eq:canonical}
y^* = y + \lambda \ , \qquad w^* = w \ ,
\end{equation}
or the so called normal variables, see Fazio \cite{Fazio:1990:NVT}, that obey the relations
\begin{equation}\label{eq:normal}
y^* = y + (1 +\lambda \delta) \ , \qquad w^* = w (1 + \lambda) \ .
\end{equation}
We note that (\ref{eq:normal}) is the infinitesimal form of a scaling group.
In the case of the normal variable, taking into account the equations (\ref{eq:inf2}) and (\ref{eq:normal}), we get
\begin{align}\label{eq:system:normal}
& X \frac{\partial y}{\partial x} + U \frac{\partial y}{\partial u} = \delta y \nonumber \\[-1ex]
&\\[-1ex]
& X \frac{\partial w}{\partial x} + U \frac{\partial w}{\partial u} = w \ . \nonumber 
\end{align}
The general solution of (\ref{eq:system:normal}) can be found by integrating the characteristic equations
\begin{equation}\label{eq:char}
\frac{dx}{X(x, u)} = \frac{du}{U(x, u)} = \frac{dy}{\delta y} = \frac{dw}{w} \ .
\end{equation}
By integrating the first equation in (\ref{eq:char}) we get the so-called first group invariant $I = I(x, u)$, whereupon we have
\begin{align}\label{eq:newvar}
& y = y_0(I) \exp\left(\delta \int{\frac{dx}{X(x, u(x))}}\right) \nonumber \\[-1ex]
&\\[-1ex]
& w = w_0(I) \exp\left(\int{\frac{dx}{X(x, u(x))}}\right) \ , \nonumber 
\end{align}
where $y_0(I)$ and $w_0(I)$ are arbitrary functions of the first invariant.
Fazio \cite{Fazio:1992:NLE} gives an application of the idea developed above to the numerical length estimation of tubular flow reactors.   

\vspace{1cm}
\noindent
{\bf Acknowledgements:} This work was partially supported by GNCS of INDAM and by University of Messina.


\end{document}

%% file: TABes1.tex
\begin{table}[!hbt]
\caption{Sample numerical results for the problem (\ref{eq:ex1:fbf}) with $\e = 1.0\mbox{D}-06$. Here and in the following we use the notation $\mbox{D}-k = 10^{-k}$ for a double precision arithmetic.}
\vspace{.5cm}
\renewcommand\arraystretch{1.3}
	\centering
		\begin{tabular}{cr@{.}lr@{.}lr@{.}l}
\hline 
Method 
& \multicolumn{2}{c}%
{$x_\e$}
& \multicolumn{2}{c}%
{$u_\e(0)$}
& \multicolumn{2}{c}%
{${\ds \frac{du_\e}{dx}(0)}$} \\[1.2ex]
\hline
RK4 &  13 & 8152456 &  2 & 6660$\mbox{D}-04$ & 0 & 999734403 \\
RK6 &  13 & 8152449 &  2 & 6659$\mbox{D}-04$ & 0 & 999734408 \\
RK8 &  13 & 8152449 &  2 & 6659$\mbox{D}-04$ & 0 & 999734408 \\
\hline			
		\end{tabular}
	\label{tab:RKes1}
\end{table}

%% file: TABes1Dx.tex
\begin{table}[!hbt]
\caption{Numerical results for the problem (\ref{eq:ex1:fbf}) for  $\e = 1\mbox{D}-06$ and different values of  $\Delta x$. We used RK6.}
\vspace{.5cm}
\renewcommand\arraystretch{1.3}
	\centering
		\begin{tabular}{lr@{.}lr@{.}lr@{.}l}
\hline 
{$-\Delta x$}
& \multicolumn{2}{c}%
{$x_\e$}
& \multicolumn{2}{c}%
{$u_\e(0)$}
& \multicolumn{2}{c}%
{${\ds \frac{du_\e}{dx}(0)}$} \\[1.2ex]
\hline
0.1        &  13 & 8149 &     6 & 48$\mbox{D}-04$ & 0 & 999353 \\
0.05      &  13 & 8152 &     2 & 67$\mbox{D}-04$ & 0 & 999734 \\
0.025    &  13 & 8154 &     7 & 37$\mbox{D}-05$ & 0 & 999927 \\
0.0125  &  13 & 8155 &     1 & 42$\mbox{D}-05$ & 0 & 999987 \\
0.00625  &  13 & 8155 &    4 & 88$\mbox{D}-06$ & 0 & 999996 \\
0.003125  &  13 & 8155 &   1 & 71$\mbox{D}-07$ & 1 & 000001 \\
\hline			
		\end{tabular}
	\label{tab:RKes1:Dx}
\end{table}

%% file: TABes1eps.tex
\begin{table}[!hbt]
\caption{Numerical results for the problem (\ref{eq:ex1:fbf}) for  $\Delta x = -1.0\mbox{D}-03$ and different values of  $\e$. We used RK6.}
\vspace{.5cm}
\renewcommand\arraystretch{1.3}
	\centering
		\begin{tabular}{lr@{.}lr@{.}lr@{.}l}
\hline 
{$\e$}
& \multicolumn{2}{c}%
{$x_\e$}
& \multicolumn{2}{c}%
{$u_\e(0)$}
& \multicolumn{2}{c}%
{${\ds \frac{du_\e}{dx}(0)}$} \\[1.2ex]
\hline
1.0$\mbox{D}-01$ &    2 & 39790 &     5 & 16$\mbox{D}-08$ & 1 & 099999948 \\
1.0$\mbox{D}-02$ &    4 & 61512 &     5 & 35$\mbox{D}-08$ & 1 & 009999946 \\
1.0$\mbox{D}-03$ &    6 & 90875 &     9 & 26$\mbox{D}-08$ & 1 & 000999907 \\
1.0$\mbox{D}-04$ &    9 & 21044 &     1 & 23$\mbox{D}-07$ & 1 & 000099877 \\
1.0$\mbox{D}-05$ &  11 & 5129 &    3 & 02$\mbox{D}-08$ & 1 & 000009970 \\
1.0$\mbox{D}-06$ &  13 & 8155 &   1 & 25$\mbox{D}-07$ & 1 & 000000875 \\
1.0$\mbox{D}-07$ &  16 & 1181 &   4 & 33$\mbox{D}-08$ & 1 & 000000057 \\
1.0$\mbox{D}-08$ &  18 & 4207 &   1 & 09$\mbox{D}-07$ & 0 & 999999901 \\
1.0$\mbox{D}-09$ &  20 & 7233 &   9 & 76$\mbox{D}-08$ & 0 & 999999903 \\
\hline			
		\end{tabular}
	\label{tab:RKes1:eps}
\end{table}

%% file: TABes2.tex
\begin{table}[!hbt]
\caption{Sample numerical results for the problem (\ref{eq:ex2:fbf}).}
\vspace{.5cm}
\renewcommand\arraystretch{1.3}
	\centering
		\begin{tabular}{cr@{.}lr@{.}lr@{.}l}
\hline 
Method &
\multicolumn{2}{c}%
{${\ds \frac{du_\e}{dx}(x_\e)}$}
& \multicolumn{2}{c}%
{$x_\e$}
& \multicolumn{2}{c}%
{${\ds \frac{du_\e}{dx}(0)}$} \\[1.2ex]
\hline
RK4 &  8 & 24683$\mbox{D}-09$ &      9 & 999999539 & 1 & 000000092 \\
RK6 &  8 & 24462$\mbox{D}-09$ &     10 & 000000058 & 0 & 999999988 \\
RK8 &  8 & 24461$\mbox{D}-09$ &     9  & 999999959 & 1 & 000000008 \\
\hline			
		\end{tabular}
	\label{tab:RKes2}
\end{table}